\documentclass[12pt]{article}

\usepackage{amsmath,amssymb,amsfonts,theorem,makeidx,latexsym,epsfig,subfigure}

\usepackage{color}

\newtheorem{defn}{Definition}[section]

\newtheorem{lemma}[defn]{Lemma}

{\theorembodyfont{\rmfamily}

\newtheorem{ex}[defn]{Example}}

\newtheorem{thm}[defn]{Theorem}

\newtheorem{prop}[defn]{Proposition}

\newtheorem{cor}[defn]{Corollary}

\newtheorem{rem}[defn]{Remark}

\numberwithin{equation}{section}

\newcommand{\h}{{\cal H}}

\newcommand{\ltr}{ L^2(\mathbb R) }

\newcommand{\mn}{\mathbb N}

\newcommand{\mr}{\mathbb R}

\newcommand{\mz}{\mathbb Z}

\def\bp{{\noindent\bf Proof. \ }}

\def\ep{\hfill$\square$\par\bigskip}

\def\bqs{\begin{equation}}

\def\eqs{\tag*{$\square$}\end{equation}\par\bigskip}

\def\hp{\hat{\psi}}

\def\la{\langle}

\def\ra{\rangle}

\def\ga{\gamma}

\def\sukz{\sum_{k\in \mz}}

\def\sujz{\sum_{j\in \mz}}

\def\hp{\hat{\psi}}

\def\nl{\left|\left|}

\def\nr{\right|\right|}

\def\bop{\begin{op}\rm}

\def\eop{\end{op}}

\def\bee{\begin{eqnarray}}

\def\ene{\end{eqnarray}}

\def\bes{\begin{eqnarray*}}

\def\ens{\end{eqnarray*}}

\def\bei{\begin{itemize}}

\def\eni{\end{itemize}}

\def\bt{\begin{thm}}

\def\et{\end{thm}}

\def\bc{\begin{cor}}

\def\ec{\end{cor}}

\def\bpr{\begin{prop}}

\def\epr{\end{prop}}

\def\bl{\begin{lemma}}

\def\el{\end{lemma}}

\def\bd{\begin{defn}}

\def\ed{\end{defn}}

\def\bex{\begin{ex}}

\def\enx{\end{ex}}

\def\bfi{\begin{fig}}

\def\efi{\end{fig}}

\newcommand{\nf}{\nl f \nr}

\def\hp{\widehat{\psi}}

\def\ftki{\{f_k\}_{k\in I}}
\def\gtki{\{g_k\}_{k\in I}}
\def\suki{\sum_{k\in I}}

\title{Approximately dual pairs of wavelet frames}

\date{\today}

\author{Ana Benavente, Ole Christensen, Marzieh Hasannasab \\ Hong Oh Kim, Rae Young Kim,
Federico D. Kovac}

\begin{document}

\maketitle

\begin{abstract} This paper deals with structural issues concerning wavelet frames and their dual frames. It is known that there exist wavelet frames
$\{a^{j/2}\psi( a^j\cdot -kb)\}_{j,k\in \mz}$ in $\ltr$ for which no dual frame
has wavelet structure. We first generalize this result by proving that there exist wavelet
frames for which no approximately dual frame has wavelet structure. Motivated
by this
we show that by imposing a very mild decay condition on
the Fourier transform of the generator $\psi\in \ltr,$
a certain oversampling $\{a^{j/2}\psi( a^j\cdot -kb/N)\}_{j,k\in \mz}$  indeed has an approximately dual
wavelet frame; most importantly, by choosing the parameter $N\in \mn$ sufficiently
large we can get as close to perfect reconstruction as desired, which makes the approximate dual frame pairs perform
equally well as the classical dual frame pairs in applications.
\end{abstract}

\begin{minipage}{120mm}


\end{minipage}
\

\section{Introduction}
Given any $a>1$ and $b\in \mr,$ consider the associated {\it scaling operator}
$D_a$ and translation operator $T_b$ acting on $\ltr$ by $D_af(x):=a^{1/2}f(ax),$
resp., $T_bf(x)=f(x-b), \, x\in \mr.$ The {\it wavelet system} generated by a
function $\psi \in \ltr$ is the collection of functions
$\{D_{a^j} T_{kb}\psi\}_{j,k\in \mz}=\{a^{j/2}\psi( a^j\cdot -kb)\}_{j,k\in \mz}.$
It is well-known \cite{Da1,CS} that even though a wavelet frame always
has a dual frame, there are cases where no dual frame with wavelet structure  exists. Thus,
there might not exist a function $\widetilde{\psi}\in \ltr$ such that the
{\it perfect reconstruction formula}
\bee \label{190620a} f= \sum_{j,k\in \mz} \la f, D_{a^j} T_{kb}\widetilde{\psi}\ra D_{a^j} T_{kb}\psi \ene holds for all $f\in \ltr.$ For this reason we will
consider a weakening of the condition \eqref{190620a}; indeed, following \cite{CLau},
we say that two wavelet systems $\{D_{a^j} T_{kb}\psi\}_{j,k\in \mz}$ and
$\{D_{a^j} T_{kb}\widetilde{\psi}\}_{j,k\in \mz}$ form {\it approximately
dual frames} if both systems form Bessel sequences, and there exists $\epsilon<1$ such that
\bee \label{190620b} \nl  f-\sum_{j,k\in \mz} \la f, D_{a^j} T_{kb}\widetilde{\psi}\ra D_{a^j} T_{kb}\psi \nr  \le \epsilon\, ||f||, \, \forall f\in \ltr. \ene
Since real-life applications always involve some truncations and approximations,
approximately dual frames
with sufficiently small parameters $\epsilon$
( say, smaller
than machine precision)  will perform as well as the classical dual pairs of frames. At
the same time it is a significantly weaker condition to ask for
functions $\psi, \widetilde{\psi}$ to generate approximately dual frames, so
giving up the perfect reconstruction property leads to a significant gain
in design freedom.
In our main result
we prove that under a very mild decay condition
on the Fourier transform of a function
$\psi \in \ltr$ that generates a wavelet
frame $\{D_{a^j}T_{kb}\psi\}_{j,k\in \mz},$ there
exists $N\in \mn$ such that the  {\it oversampled wavelet system} $\{D_{a^j}T_{kb/N}\psi\}_{j,k\in \mz}$
has an approximately dual wavelet frame; most importantly, by
choosing $N$ sufficiently large, we can get as
close to perfect reconstruction as desired.

Approximate duals already appear in a number
of contexts in the literature, see
the papers \cite{BuiLau,BCZ,CKK19,DMa2,FGry,FOnc}.
Also,  the parallel case of Gabor frames
was considered in the paper \cite{CKKJ}. Even though
some of the techniques are similar, there is a fundamental
difference between the Gabor case and the wavelet
case considered in the current paper: indeed,
in contrast to the wavelet case, a Gabor
frame always has an exact dual frames with the same
structure.  For the wavelet case, there is an
extensive literature concerned with the case of
(multi)wavelets where dual frames with the same
structure exists, typically formulated within
the setting of the unitary extension principle
and its later variants.
However, these constructions are based on refinable functions and thus yield a quite limited class.
General dual wavelet frames have been considered
and characterized, e.g., in \cite{FGWW,Bow2,HLW},
but the literature only contains sparse information
about the case where a wavelet frames does not
have a dual with the same structure.

Our first result, stated in Section \ref{190620d}, is negative and will serve
as motivation for the rest of the paper: we prove
that there exist wavelet frames, where not even an  approximately dual wavelet system exist,
i.e., giving up the perfect reconstruction property by itself does not lead to
sufficient extra design freedom.
Our main new contribution, outlined above, is then developed in Section \ref{190620e}

\section{On the existence of approximately dual wavelet frames} \label{190620d}

The purpose in this Section is to show that there exist wavelet frames for
which no approximate dual wavelet frame exists. The result will be derived
as a consequence of a couple of results  that are derived in the setting of frames for a general Hilbert space $\h.$ Since we are aiming at a no-go result, it is enough
to consider a special class of frames, and our results will actually be derived
in the setting of Riesz bases. The first result yields a lower bound for
the deviation from perfect reconstruction that is obtained when the dual Riesz
basis is replaced by any other Bessel sequence in the Hilbert space. For the
sake of our application to wavelet systems we will index the Riesz basis by
an arbitrary countable index set $I.$

\bpr \label{70406a} Let $\ftki$ be a Riesz basis for a Hilbert space $\h,$ with  dual Riesz basis $\gtki$ and lower frame
bound $A.$ Then, for any Bessel sequence $\{ \varphi_k\}_{k\in I}$ and
any $\ell\in \mn,$
\bee \label{92408a} \sup_{||f||=1} \nl f- \sum_{k\in I} \la f, \varphi_k\ra f_k \nr^2 \ge A \, ||g_\ell- \varphi_\ell||^2.\ene \epr
\bp For any $f\in \h,$ using that $f= \suki \la f, g_k\ra f_k,$ we obtain that

\bes \nl f- \suki \la f, \varphi_k\ra f_k \nr^2  =
\nl \suki \la f, g_k- \varphi_k\ra f_k \nr^2  \ge   A
\suki | \la f, g_k- \varphi_k\ra |^2 .\ens Now, if
$g_\ell=\varphi_\ell,$ the result \eqref{92408a} clearly holds;
otherwise,   taking
$f:= (g_\ell - \varphi_\ell)/||g_\ell - \varphi_\ell||$
yields the result. \ep

As a consequence of
Proposition \ref{70406a} we now show that if   the vectors $\gtki$ in the dual Riesz basis do not
have equal norm, then a sequence
$\{ \varphi_k\}_{k\in I}$ of vectors with identical norms can not be an
approximate dual frame for arbitrary small values of $\epsilon:$

\bc\label{70406b} Let $\ftki$ be a Riesz basis with  dual Riesz basis $\gtki$ and lower frame
bound $A.$ Assume that there exists $j,\ell\in I$ such that $\delta:= ||g_j||- ||g_\ell|| >0,$ and consider any Bessel sequence $\{ \varphi_k\}_{k\in I}$ consisting of vectors with equal norm. Then
\bes \sup_{||f||=1} \nl f- \suki \la f, \varphi_k\ra f_k \nr^2 \ge \frac{A\delta^2}{4}.\ens \ec

\bp Choose $j,\ell\in I$ such that $\delta:= ||g_j||- ||g_\ell|| >0.$ If the vectors $\{ \varphi_k\}_{k\in I}$ have equal norm, then  either
$|| g_j- \varphi_j||\ge \delta/2$ or $|| g_\ell- \varphi_\ell||\ge \delta/2$. Thus
the result follows directly from Proposition \ref{70406a}. \ep

The relevance of Corollary \ref{70406b} for the aim of the section is
immediately clear: indeed, if we take a Riesz basis $\ftki$ with
wavelet structure and we aim at another wavelet system $\{ \varphi_k\}_{k\in I}$
being an approximately dual frame, then all the vectors in $\{ \varphi_k\}_{k\in I}$
have equal norm and we can apply the result in Corollary \ref{70406b} directly.
Using this idea we will now exhibit a class of wavelet frames, where
it is impossible to get arbitrarily close to perfect reconstruction
using another wavelet system.
The frame construction itself  originally appeared in \cite{Da1,CS} and was
elaborated further in \cite{CB}.
\bex \label{190620g} Let $\{D_{2^j} T_{k}\psi\}_{j,k\in \mz}$ be a wavelet
orthonormal basis for $\ltr$. Given $\eta \in ]0,1[$,
define a function $\theta$ by \bes \theta = \psi + \eta D_2 \psi.
\ens In \cite{Da1} it is proved that  $\{D_{2^j}T_k \theta\}_{j,k\in \mz}$ is a
Riesz basis and that the dual Riesz basis does not have wavelet structure.
Using that the scaling operator $D_2$ is unitary and that $\{D_{2^j} T_{k}\psi\}_{j,k\in \mz}$ is an orthonormal system, it follows that for any finite
scalar sequence $\{c_{j,k}\}_{(j,k)\in I}, \, I \subset \mz \times \mz,$
\bes \nl \sum_{(j,k)\in I} c_{j,k} \left( D_{2^j}T_k\theta - D_{2^j}T_k \psi\right) \nr & = & \nl \sum_{(j,k)\in I} c_{j,k} D_{2^j}T_k (\eta D_2 \psi) \nr \\
& = & \eta \, \nl \sum_{(j,k)\in I} c_{j,k} D_{2^{j+1}}T_{2k}  \psi \nr =
\eta \left( \sum_{(j,k)\in I} |c_{j,k}|^2\right)^{1/2}.
\ens
Using standard perturbation results for frames (see, e.g., Theorem 22.1.1 in \cite{CB}),
it follows that the Riesz basis $\{D_{2^j}T_k \theta\}_{j,k\in \mz}$
has lower frame bound $A=(1- \eta)^2.$ Let $\{\omega_{j,k}\}_{j,k\in \mz}$
denote
the dual Riesz basis associated with $\{D_{2^j}T_k \theta\}_{j,k\in \mz}$.
In Example 16.1.1 in \cite{CB}
it is proved that for any $j,k\in \mz$ we have $\omega_{j, 2k+1}=D_{2^j}T_{2k+1}\psi;$
thus $|| \omega_{j, 2k+1}||=1 $ for $j,k\in \mz.$ It is also proved that
$\omega_{j,0}= \sum_{n=0}^\infty (-\eta)^n D_{2^{j-n}}\psi$ for any $j\in \mz,$
which implies that $||\omega_{j,0}||= (1- \eta^2)^{-1/2}.$ In particular, we see that
$\delta:= ||\omega_{0,0}|| - ||  \omega_{0,1}||= (1- \eta^2)^{-1/2}-1.$ Now,
consider any wavelet system $\{D_{2^j}T_k \varphi\}_{j,k\in \mz}$ in $\ltr.$
Since a wavelet system always consists of
vectors of equal norm, Corollary \ref{70406b} immediately implies that

\bes \sup_{||f||=1} \nl f - \sum_{j,k\in \mz} \la f, D_{2^j}T_k \varphi\ra D_{2^j} T_k \theta\nr^2
\ge \frac{(1- \eta)^2}{4} \left( \frac1{\sqrt{1-\eta^2}} -1 \right)^2.\ens
This estimate immediately shows that for the Riesz basis  $\{D_{2^j}T_k \theta\}_{j,k\in \mz}$ there is a limition on how close we can get to perfect reconstruction using
an ``approximate dual wavelet system" $\{D_{2^j}T_k \varphi\}_{j,k\in \mz}$ in $\ltr.$
\ep \enx

\section{Construction of approximately dual frames} \label{190620e}

Motivated by the negative result in Example \ref{190620g} we will now
develop a new strategy for obtaining pairs of approximately dual wavelet frames.
The key idea is to use {\it oversampling.} Indeed, our main result
shows that for a large class of wavelet systems
$\{D_{a^j} T_{kb}\psi\}_{j,k\in \mz}$ we can get arbitrarily
close to perfect reconstruction by using two wavelet systems that are obtained from
the given one using oversampling.

 A key step in our construction of an approximately dual wavelet frames
is to truncate the given function $\psi \in \ltr.$ Throughout the paper,
given $\psi\in \ltr$ and any $K>0,$ define the function
$\psi_K \in \ltr$ by
\bee\label{190222b} \widehat{\psi_K}:= \hp \, \chi_{[-K,K]}.\ene

We will need the following result, which basically estimates the Bessel bound of the wavelet
system $\{D_{a^j}T_{kb}(\psi- \psi_K)\}_{j,k\in \mz}.$

\newpage

\bpr\label{190612d} Let $a>1$ and $b>0$. Let $\psi \in
\ltr$ and assume that there exist constants $C>0$ and $\sigma>0$
such that \bee \label{191017d} |\hp(\ga)|\le  \frac{C}{1+ |\ga |^{1+\sigma}}, \ \ \, \ga \in \mr. \ene
Then
\begin{equation*}
\sup_{|\ga| \in [1,a]}  \sum_{j,k\in \mz}
|(\hp-\widehat{\psi_K})(a^j \ga) (\hp-\widehat{\psi_K})(a^j \ga +
k/b)|
 \le \frac{2C^2}{\sigma K^{1+2\sigma}}\left( \frac{\sigma}{K}+b\right)
\frac{a^{1+\sigma}}{a^{1+\sigma}-1}.
\end{equation*}
\epr

\bp   For notational convenience, let $p_K(x) := \psi(x)-
\psi_K(x).$ We first note that
\begin{eqnarray}
&&\sup_{|\ga| \in [1,a]}  \sum_{j,k\in \mz}
|\widehat{p_K}(a^j \ga) \widehat{p_K}(a^j \ga +
k/b)| \nonumber \\
&& \le
\left(\sup_{|\ga| \in [1,a], j\in \mz}
\sum_{k \in \mz}
\bigg|\widehat{p_K}(a^j
\ga+k/b)\bigg|\right)
\left(\sup_{|\ga| \in [1,a]}
\sum_{j\in \mz} \bigg| \widehat{p_K}(a^j\ga)\bigg|\right) \label{190612a}
\end{eqnarray}
and $\widehat p_K(\ga)
=  \widehat \psi (\ga) \chi_{[-K,K]^c}(\ga), \, \ga\in \mr.$ We now proceed
with a number of estimates that finally yields the result.

\vspace{.1in}\noindent
{\bf Estimate of $\sum_{k \in \mz} \bigg| \widehat p_N(a^j \ga+k/b)\bigg|$ :}
Fix $\ga\in[-a,-1]\cup[1,a]$ and $j   \in\mz$. Then
the contribution from
$a^j \ga+k/b$ hitting the interval $]-\infty,-K]\cup[K,\infty[$ is at most
\begin{eqnarray*}
&& \sum_{k=0}^{\infty }
\left|\hp\left(K+\frac{k}{b}\right)\right|+
\sum_{k=0}^{\infty } \left|\hp\left(-K-\frac{k}{b}\right)\right|\\
&&\leq  2\sum_{k=0}^{\infty }
\frac{C}{\left(K+\frac{k}{b}\right)^{1+\sigma}} \nonumber \\
&&\leq  \frac{2C}{K^{1+\sigma}}+ \int_0^\infty
\frac{2C}{\left(K+\frac{t}{b}\right)^{1+\sigma}}dt \\
&&=2C \left(\frac{1}{K^{1+\sigma}}+\frac{b}{\sigma K^{\sigma}} \right) \\
&&=\frac{2C}{\sigma K^{\sigma}}\left( \frac{\sigma}{K}+b\right).
\end{eqnarray*}
This leads to
\begin{equation}\label{190612b}
\sum_{k \in \mz} \bigg|  \widehat p_{N}(a^j \ga+k/b) \bigg|\leq
 \frac{2C}{\sigma K^{\sigma}}\left( \frac{\sigma}{K}+b\right).
\end{equation}

\noindent {\bf Estimate of $\sum_{j\in \mz} \bigg|  \widehat p_N(a^j\ga)\bigg|$ :}
Fix $\ga\in[-a,-1]\cup [1,a]$. Then  the contribution from
$a^j \ga$ hitting the interval $]\infty,-K]\cup[K,\infty[$ is at most
\begin{eqnarray*}
&&\max\left\{ \sum_{j=0}^{\infty } \left| \hp\left(a^j
K\right)\right|,
 \sum_{j=0}^{\infty } \left| \hp\left(-a^j K\right)\right|\right\}\\
 &&\leq  \sum_{j=0}^{\infty }
\frac{C}{\left(a^j K\right)^{1+\sigma}} = \frac{C}{K^{1+\sigma}}\sum_{j=0}^\infty
\frac{1}{(a^{1+\sigma})^j}
=\frac{C}{K^{1+\sigma}}\frac{1}{1-\frac{1}{a^{1+\sigma}}}
=\frac{Ca^{1+\sigma}}{K^{1+\sigma}\left(a^{1+\sigma}-1\right)}.
\end{eqnarray*}
This leads to
\begin{equation}\label{190612c}
\sum_{j\in \mz} \bigg|  \widehat p_N(a^j\ga)\bigg|\leq
\frac{Ca^{1+\sigma}}{K^{1+\sigma}\left(a^{1+\sigma}-1\right)}.
\end{equation}
The  result now follows from \eqref{190612a} together with \eqref{190612b} and
\eqref{190612c}.
\ep

Proposition \ref{190612d} implies that if
$\{D_{a^j}T_{kb}\psi\}_{j,k\in \mz}$ is a
frame and \eqref{191017d} holds, then also
$\{D_{a^j}T_{kb}\psi_K\}_{j,k\in \mz}$ is a
frame whenever $K\in \mn$ is sufficiently large.
We are now ready to present our main result.

\bt \label{202003ag}  Assume that for some $a>1, b>0,$ the function $\psi\in \ltr$ generates a wavelet frame and
that  for some $C>0, \sigma>0,$
\bee \label{191017dn} |\hp(\ga)|\le  \frac{C}{1+ |\ga |^{1+\sigma}}, \ \ \, \ga \in \mr. \ene
Given any $K\in \mn,$ let $N$ denote the unique odd
integer such that
\bee \label{202003b} N-2 \le 2bK <N,\ene
and define
the function $\widetilde{\psi_K}$ by
 \bee\label{190222c} \widehat{\widetilde{\psi_K}}(\ga) =
\frac{b\widehat{\psi_K}(\ga)}{N\sum_{k\in \mz}| \widehat{\psi_K}(a^k
\ga)|^2}.\ene
Then, letting
\bee \label{202003c} \varepsilon_K:=
 \left(\frac{K+1/b}{K^{3+2\sigma}}\right)^{1/2}
2C \left[(1+ \frac{1}{2\sigma})\frac{a^{1+\sigma}}{a^{1+\sigma}-1}  \right]^{1/2} \ene
and choosing
$K\in \mn$ large enough such that $\varepsilon_K <\sqrt{2bA},$  we have
\bee \label{202003k} \nl f - \sum_{j,k\in \mz}
\la f, D_{a^j}T_{kb/N} \widetilde{ \psi_K}\ra
D_{a^j}T_{kb/N}\psi\nr \le\frac{ \varepsilon_K}{\sqrt{2bA}- \varepsilon_K}\, \nf, \, \forall
f\in \ltr.\ene \et

\bp
Using Proposition \ref{190612d} on the function $\psi-\psi_K$
and with the translation parameter $b$ replaced by
$b/N$, \eqref{202003b} implies that
\bes & \ & \frac{N}{b} \sukz \sujz | (\widehat{\psi}-\widehat{\psi_K})(a^j\ga)
(\widehat{\psi}-\widehat{\psi_K})(a^j\ga +kN/b)| \\
 &  &\le \frac{N}{b} \frac{2C^2}{\sigma K^{1+2\sigma}}\left( \frac{\sigma}{K}+ \frac{b}{N}\right)
\frac{a^{1+\sigma}}{a^{1+\sigma}-1} \\
 &  &=
\frac{4C^2}{K^{2+2\sigma}}\left(\frac{N}{2b}+\frac{K}{2\sigma} \right) \frac{a^{1+\sigma}}{a^{1+\sigma}-1} \\
&  &\le \frac{4C^2}{K^{2+2\sigma}}
\left(K+\frac{1}{b}+\frac{K}{2\sigma} \right)\frac{a^{1+\sigma}}{a^{1+\sigma}-1} \\
&  &\le
\frac{4C^2}{K^{2+2\sigma}}\left(K+\frac{1}{b}\right)\left(1+
\frac{1}{2\sigma}\right)\frac{a^{1+\sigma}}{a^{1+\sigma}-1} =:R_K. \ens
Thus the wavelet system $\{D_{a^j}T_{kb/N}(\psi-\psi_K)\}_{j,k\in
\mz}$ is a Bessel sequence with bound $R_K.$ Now, since the
oversampling parameter $N$ is odd, the oversampling theorem by
Chui and Shi \cite{CS3} shows that the wavelet system
$\{D_{a^j}T_{kb/N}\psi\}_{j,k\in \mz}$ is a frame with lower frame
bound $AN;$ hence by standard frame perturbation theory (see  Cor.
22.1.5 in \cite{CB}), $\{D_{a^j}T_{kb/N}\psi_K\}_{j,k\in \mz}$ is
also a frame whenever $R_K < AN,$ with lower frame bound $
AN\left( 1- \sqrt{\frac{R_K}{AN}}\right)^2.$ Invoking again the
relationship between $K$ and $N$ in \eqref{202003b}, we conclude
that in particular $\{D_{a^j}T_{kb/N}\psi_K\}_{j,k\in \mz}$ is a
frame whenever $R_K < 2bAK,$ with lower frame bound $ AN\left( 1-
\sqrt{\frac{R_K}{2bAK}}\right)^2.$ By construction,
$\widehat{\psi}_K$ is supported on an interval of length $2K;$
hence, again due to \eqref{202003b} the canonical dual frame of
$\{D_{a^j}T_{kb/N}\psi_K\}_{j,k\in \mz}$ is
$\{D_{a^j}T_{kb/N}\widetilde{\psi_K}\}_{j,k\in \mz},$ with
$\widetilde{\psi_K}$ defined as in \eqref{190222c}. Note that the
upper frame bound for
$\{D_{a^j}T_{kb/N}\widetilde{\psi_K}\}_{j,k\in \mz}$ is $
\left[AN\left( 1- \sqrt{\frac{R_K}{2bAK}}\right)^2\right]^{-1}.$
Thus, for any $f\in \ltr$ we arrive at \bes & \ & \nl f -
\sum_{j,k\in \mz} \la f, D_{a^j}T_{kb/N} \widetilde{ \psi_K}\ra
D_{a^j}T_{kb/N}\psi_K\nr^2 \\
&  &= \nl \sum_{j,k\in \mz}
\la f, D_{a^j}T_{kb/N} \widetilde{ \psi_K}\ra
D_{a^j}T_{kb/N}\psi_K - \sum_{j,k\in \mz}
\la f, D_{a^j}T_{kb/N} \widetilde{ \psi_K}\ra
D_{a^j}T_{kb/N}\psi\nr^2 \\
&  &= \nl \sum_{j,k\in \mz}
\la f, D_{a^j}T_{kb/N} \widetilde{ \psi_K}\ra
D_{a^j}T_{kb/N}\left(\psi- \psi_K\right) \nr^2 \\
&  &\le R_K \sum_{j,k\in \mz}
|\la f, D_{a^j}T_{kb/N} \widetilde{ \psi_K}|^2 \\
&  &= \frac{R_K}{AN\left( 1- \sqrt{\frac{R_K}{2bAK}}\right)^2}\, ||f||^2.\ens This implies that
\bes \nl f - \sum_{j,k\in \mz}
\la f, D_{a^j}T_{kb/N} \widetilde{ \psi_K}\ra
D_{a^j}T_{kb/N}\psi_K\nr
\le \frac{\sqrt{\frac{R_K}{AN}}}{ 1- \sqrt{\frac{R_K}{2bAK}}}\, ||f|| \le  \frac{ \sqrt{\frac{R_K}{K}}}{\sqrt{2bA}-\sqrt{\frac{R_K}{K}}}\, \nf;\ens the proof is now completed
by reformulating the result in terms
of $\varepsilon_K= \sqrt{R_K/K}.$
 \ep

\begin{rem} \rm{ Theorem \ref{202003ag} yields a
very explicit relationship between the support-size
of the truncated window, measured by the parameter $K,$
and the degree of approximation; in fact,
asymptotically as $K\to\infty$ and for some constant
$D>0,$
\bee \label{203003a} \nl f - \sum_{j,k\in \mz}
\la f, D_{a^j}T_{kb/N} \widetilde{ \psi_K}\ra
D_{a^j}T_{kb/N}\psi\nr \le  \frac{D}{K^{1+\sigma}}\, ||f||, \, \forall f\in \ltr.\ene}
\end{rem}

In particular, our method allows to us get as close
as desired to perfect reconstruction using
oversampling of the given wavelet system.

\bc \label{190616f}  Under the assumptions in
Theorem \ref{202003ag}, given any $\epsilon>0,$
choose $K\in \mn$ such that
\bes\frac{ \varepsilon_K}{\sqrt{2bA}- \varepsilon_K} \le\epsilon.\ens
Then, denoting
the synthesis operator for the oversampled frame
$\{D_{a^j}T_{kb/N}\psi\}_{j,k\in \mz}$
bu $U_N$ and the synthesis operator for
$\{D_{a^j}T_{kb/N}\widetilde{\psi_K}\}_{j,k\in \mz}$
by
$W_N,$ we have that
\bee \label{190828af} || I- U_NW_N^*||\leq
\epsilon.  \ene
\ec

Let us now return to the motivating Example \ref{190620g}, and show that
the condition in Theorem \ref{202003ag}
is satisfied when we apply
the construction to the Battle-Lemari\'{e} wavelets:

\bex Let us first recall the key steps in the construction
of the Battle-Lemari\'{e} wavelets \cite{Battle, Lemarie}; see, e.g., \cite{Da2}
for more details about the construction. The B-splines $N_m, \, m\in \mn,$ are defined
recursively by $ N_1:= \chi_{[0, 1]}$ and the convolution
formula $ N_{m+1}=N_m *N_1.$ It is well known that
$\widehat{N_m}(\ga)  =
\left(\frac{1- e^{-2\pi i \ga}}{2\pi i \ga}\right)^m$ and  that $\widehat{N_m}(2\ga)=H_0(\ga) \widehat{N_m}(\ga),$ where
$ H_0(\ga)= \left(\frac{1+e^{-2\pi i \ga}}{2}\right)^m;$
furthermore,
there exist constants $A,B>0$ such that $A\le G(\ga):=\sum_{k\in
\mz} | \widehat{N_m}(\ga+k)|^2 \le B, \ \ga\in \mr.$
Define the function $\varphi\in \ltr$
by its Fourier transform $\widehat{\varphi}$ via
$\widehat{\varphi}(\ga):= \frac1{\sqrt{G(\ga)}} \,
\widehat{N_m}(\ga),\, \ga \in \mr.$
Then $\varphi$ is a scaling function; indeed, by direct calculation,
$ \widehat{\varphi}(2\ga)
= \widetilde{H_0}(\ga)\widehat{\varphi}(\ga),$ where
$ \widetilde{H_0}(\ga):= \sqrt{\frac{G(\ga)}{G(2\ga)}}\, H_0(\ga).$
The Battle-Lemari\'{e} wavelet $\psi$ associated with the scaling function $\varphi$ is given in the
Fourier domain by the formula
\bes \widehat{\psi}(\ga)  =  \overline{\widetilde{H_0}(\ga/2 + 1/2)}e^{-\pi i \ga}
\widehat{\varphi}(\ga/2)  =  \sqrt{\frac{G(\ga/2 + 1/2)}{G(\ga) G(\ga/2)}}\,
\left( \frac{1- \cos \pi \ga }{\pi i \ga} \right)^m e^{-\pi i \ga}.
\ens
By construction $\{D_{2^j}T_k \psi\}_{j,k\in \mz}$ is an orthonormal basis
for $\ltr.$  Given $\eta \in ]0,1[$,
define a function $\theta$ by $\theta = \psi + \eta D_2 \psi.$
Then, by Example \ref{190620g},  $\{D_{2^j}T_k \theta\}_{j,k\in \mz}$ is a
Riesz basis and  the dual Riesz basis does not have wavelet structure. Note
that
\bee \label{191017b} \widehat{\theta}(\ga)= \widehat{\psi}(\ga) + \frac{\eta}{\sqrt{2}}\widehat{\psi}(\ga/2). \ene
It is clear from \eqref{191017b} and the formula
for $\widehat{\psi}(\ga)$  that the function $\theta$
satisfies the decay condition in Theorem \ref{202003ag}
for $m\ge 2;$ thus, we can
construct approximately dual wavelet  frames based on
the Battle-Lemari\'{e} wavelets for all
 $m\ge 2$ and any $\eta \in ]0,1[,$ as claimed. \ep \enx

As already mentioned, the decay condition  in Theorem \ref{202003ag}  is very mild; it is much weaker than
the standard decay conditions used to guarantee
the Bessel condition for a wavelet system, see,
e.g., Theorem 15.2.3 and Lemma 15.2.5 in \cite{CB}.
However, if we remove this condition, we
 can not  be sure that the truncation $\psi_K$
 of the generator $\psi$ for a wavelet frame generates
 a wavelet frame itself, regardless
 how large $K$ is chosen. Based on  the subsequent Lemma \ref{190224da} we will provide a concrete
 example of this in Example \ref{190828a}.

\bl \label{190224da} Let $a>1$ and $b>0$. Given a function $\psi\in L^2(\mr)$ and $K>0$,
consider the function  $\psi_K$ in  \eqref{190222b}, and assume  that
$\{D_{a^j}T_{kb/N}\psi_K\}_{j,k\in \mz}$ is a frame
for some $N\in \mn.$ Then there
exists $J_0\in \mz$ such that
\bee \label{ff}\inf_{| \gamma | \in [1,a]}
\sum_{j=-\infty}^{J_0} \big|\widehat{\psi}(a^{j}{\gamma})\big|^{2}
>0.\ene
\el
\bp
Choose $J_0$ with $K <a^{J_0+1}$.
Due to the compact support of $\widehat{\psi_K}$, we have
$$\inf_{| \gamma | \in [1,a]} \sum_{j\in\mz} \big|\widehat{\psi_K}(a^{j}{\gamma})\big|^{2}
=\inf_{| \gamma | \in [1,a]}
\sum_{j=-\infty}^{J_0} \big|\widehat{\psi_K}(a^{j}{\gamma})\big|^{2}
\leq \inf_{| \gamma | \in [1,a]}
\sum_{j=-\infty}^{J_0} \big|\widehat{\psi}(a^{j}{\gamma})\big|^{2}.$$
Now, since $\{D_{a^j}T_{kb/N}\psi_K\}_{j,k\in \mz}$
is a frame, a result by Chui and Shi \cite{CS} says that
$ \inf_{|\ga| \in [1,a]} \sum_{j\in \mz}  \big|\widehat{\psi}(a^{j}{\gamma})\big|^{2} >0.;$
hence the desired conclusion follows.
\ep

\begin{ex} \label{190828a}
Let $a=2, b=1,$ and let $\{A_n\}_{n=0}^\infty$ be a sequence of measurable sets
such that $\left[\frac{1}{2},1\right[=\biguplus_{n=0}^\infty A_n $, where
$\biguplus$ denotes a disjoint union.
Note that  for $n\geq 1$, we have that
$2^{-n}A_n \subset \left[0,\frac{1}{2}\right[$.
Following \cite[Example 4.5(xi)]{DL}, consider the sets
\bes B:=\biguplus_{n=0}^\infty 2^{-n}A_n, \, \, C:=\left[-1, -\frac{1}{2} \right[\setminus(\biguplus_{n=1}^\infty (2^{-n}A_n-1)), \, \, D:=\biguplus_{n=1}^\infty 2^n(2^{-n}A_n-1),\ens
and define the function $\psi\in \ltr$  by $
\widehat{\psi}
:=\chi_{B\cup C\cup D}.$

It is shown in \cite[Example 4.5(xi)]{DL} that $\{D_{2^j}T_{k}\psi\}_{j,k\in \mz}$ is an orthonormal basis. Now, if  $J_0< -1$, then
\begin{equation*}
\sum_{j=-\infty}^{J_0} \big|\widehat{\psi}(2^{j}{\gamma})\big|^{2}
=0,\ \forall \ga\in [-2,-1[;
\end{equation*}
and if $J_0\geq -1$, then
\begin{equation*}
\sum_{j=-\infty}^{J_0} \big|\widehat{\psi}(2^{j}{\gamma})\big|^{2}
=0,\ \forall \ga\in \cup_{n=J_0+1}^\infty \left(2^{-n}A_{n+1}-2 \right)
\left(\subset
\left[-2,-1\right[\right).
\end{equation*}
Thus, regardless of the choice of $J_0\in\mz$, the condition \eqref{ff}
is violated. Via Lemma \ref{190224da} we conclude
that
$\{D_{a^j}T_{kb/N}\psi_K\}_{j,k\in \mz}$ can not be a frame, regardless how the parameters $K,N\in \mn$
are chosen.
\ep
\end{ex}

\noindent {\bf Acknowledgements:} This research was supported by
the Basic Science Research Program through the National Research
Foundation of Korea (NRF) funded by the Ministry of Education
(2016R1D1A1B02009954); and the 2018 Yeungnam University Research
Grant.


{\noindent Ana Benavente \\
Instituto de Matematica Aplicada San Luis \\
Universidad Nacional de San Luis and CONICET, \\
Av. Ejercito de Los Andes 950 \\
D5700 HHW San Luis, Argentina \\
E-mail: abenaven@unsl.edu.ar

\vspace{.1in} \noindent Ole Christensen\\
Department of Applied Mathematics and Computer Science\\
Technical University of Denmark,
Building 303,
2800 Lyngby, Denmark\\
Email: ochr@dtu.dk

\vspace{.1in}\noindent Marzieh Hasannasab \\
Institut f\"ur Mathematik\\
TU Berlin\\
Stra\ss e des 17. Juni 136, 10623 Berlin, Germany
\\
Email: hasannas@math.tu-berlin.de

\vspace{.1in}\noindent Hong Oh Kim \\
Department of Mathematical Sciences, KAIST\\
291 Daehak-ro, Yuseong-gu, Daejeon 34141,
Republic of Korea\\
Email: kimhong@kaist.edu

\vspace{.1in} \noindent Rae Young Kim \\
Department of Mathematics, Yeungnam University\\
280 Daehak-Ro, Gyeongsan, Gyeongbuk 38541,
Republic of Korea\\
Email:  rykim@ynu.ac.kr

\vspace{.1in} \noindent Federico D. Kovac \\
Facultad de Ingenieria
\\ Universidad Nacional de La Pampa
\\ Calle 110 no. 390 - (6360) General Pico - La Pampa,
Argentina \\
E-mail: kovacf@ing.unlpam.edu.ar

}

\end{document}